\documentclass{article}
\usepackage{amssymb}
\newtheorem{thm}{Theorem}
\newtheorem{lem}{Lemma}
\newtheorem{cor}{Corollary}
\begin{document}
\begin{center}
{\bf
A METHOD OF GENERATING SUMS OF LIKE POWERS}\\
\vspace{3ex}

\noindent \v Zarko Mijajlovi\' c$^1$, Milo\v s Milo\v sevi\'
c$^2$, Aleksandar Perovi\' c$^2$
\end{center}
\vspace{2ex}

{\footnotesize \noindent 1. Faculty of Mathematics, University of Belgrade, Studentski trg 16, 11000 Belgrade, Serbia

\hspace{-2.2ex} and Montenegro, zarkom@eunet.yu

\noindent 2. Mathematical Institute SANU, Kneza Mihaila 35, 11000 Belgrade, Serbia and Montenegro,

\hspace{-2.2ex} mionamil@eunet.yu, peramail314@yahoo.com }
\begin{abstract}
In this article a new method of generating sums of like powers is
presented.
\end{abstract}
\vspace{3ex}

We use symbols $\mathbb{Z}$, $\mathbb{N}$, $\mathbb{N}^+$,
$\mathbb{R}$ and $\mathbb{R}^+$ to represent the sets of integers,
nonnegative integers, positive integers, real numbers and positive
real numbers, respectively . In addition, we also adopt the
convention that $0^0=1$. \vspace{1ex}

Finite disjoint subsets $U$ and $V$ of $\mathbb{Z}$ will be called
a Prouhet-Tarry-Escott pair for the given integer $n>1$ if they
have the same cardinality and
$$
\sum_{u\in U}u^s=\sum_{v\in V}v^s,\ \ s=0,\dots,n-1,\ \ {\rm and}\ \
\sum_{u\in U}u^n\not=\sum_{v\in V}v^n.
$$

The notion of a $P$-sequence is recursively defined as follows:
\begin{itemize}
\item $\langle 1,-1\rangle$ is a $P$-sequence;
\item If $\langle a_0,\dots, a_k\rangle$ is a $P$-sequence and $a_0=-a_k$, then
$$
\langle a_0,\dots,a_k,a_k,\dots,a_0\rangle
$$
is also a $P$-sequence;
\item If $\langle a_0,\dots, a_k\rangle$ is a $P$-sequence and $a_0=a_k$, then
$$
\langle a_0,\dots,a_{k-1},0,-a_{k-1},\dots,-a_0\rangle
$$
is also a $P$-sequence;
\item Each $P$-sequence can be obtained only by finite use of the above clauses.
\end{itemize}
We denote the $n$-th $P$-sequence by $P_n$ (assuming that they are
ordered by their increasing lengths). For instance,
$$
P_1=\langle 1,-1\rangle,\  P_2=\langle 1,-1,-1,1\rangle,\
P_3=\langle 1,-1,-1,0,1,1,-1\rangle\ {\rm etc}.
$$
For an arbitrary positive integer $n$, $n$-th $P$-sequence $P_n
=\langle a_0,\dots,a_k\rangle$ and any integer $s\geqslant 0$ let us define a polynomial function
$F_{n,s}(x)$ by
$$
F_{n,s}(x)=\sum_{i=0}^ka_i(i+x)^s.
$$
\begin{lem}
$F_{n,s}\equiv 0$ for $s=0,\dots,n-1$.
\end{lem}
{\bf Proof.}
Since
$$
F_{n,s}(x)=((n-1)\cdots (n-1-s))^{-1}F_{n,n-1}^{(n-1-s)}(x)
$$
(here $F_{n,n-1}^{(n-1-s)}$ is the $(n-1-s)$-th derivative of $F_{n,n-1}$),
it is sufficient to prove that
\begin{eqnarray}
F_{n,n-1}\equiv 0.
\end{eqnarray}
We use induction on $n$. Trivially (1) is true for $n=1$, so let
us assume that for some $n\geqslant 1$ equality (1) holds . We
have the following two cases:
\begin{itemize}
\item $n=2m$. Assuming that $P_{2m}=\langle a_0,\dots,a_k\rangle$, we have that
$$
P_{2m+1}=\langle a_0,\dots,a_{k-1},0,-a_{k-1},\dots,-a_0\rangle
$$
and
\begin{eqnarray*}
F_{2m+1,2m}(x)&=&\sum_{i=0}^{k-1}a_i(i+x)^{2m}-\sum_{i=0}^{k-1}a_i(2k-i+x)^{2m}\\
&=&\sum_{i=0}^ka_i(i+x)^{2m}- \sum_{i=0}^ka_i(i-2k-x)^{2m}\\
&=& F_{2m,2m}(x)-F_{2m,2m}(-2k-x).
\end{eqnarray*}
Then
\begin{eqnarray*}
F_{2m+1,2m}'(x)&=& 2m\underbrace{F_{2m,2m-1}(x)}_{=0}
 - 2m\underbrace{F_{2m,2m-1}(-2k-x)}_{=0} \\
&=&0,
\end{eqnarray*}
so $F_{2m+1,2m}$ is a constant function. Since $F_{2m+1,2m}(-k)=0$, we conclude
that $F_{2m+1,2m}\equiv 0$.
\item $n=2m+1$. Similarly to the previous case one can easily check that $F_{2m+2,2m+1}$
is a constant function. Since $F_{2m+2,2m+1}(-\frac{2k+1}{2})=0$, we conclude that
$F_{2m+2,2m+1}\equiv 0$ as well. \hfill $\square$
\end{itemize}
\begin{thm}
For $s\geqslant n$ the degree of $F_{n,s}$ is equal to $s-n$.
\end{thm}
{\bf Proof.} Clearly, it is sufficient to prove that
\begin{eqnarray}
F_{n,n}\equiv const.\not=0.
\end{eqnarray}
Observe that an immediate consequence of Lemma 1 is the fact that each $F_{n,n}$ is a constant function.

We use by induction on $n$. $F_{1,1}\equiv -1$, so let us assume
that for some $n\geqslant 1$ relation (2) holds.
\begin{itemize}
\item $n=2m$. Assuming that $P_{2m}=\langle a_0,\dots, a_k\rangle$, we have that
$$
P_{2m+1}=\langle a_0,\dots,a_{k-1},0,-a_{k-1},\dots,-a_0\rangle
$$
and
\begin{eqnarray*}
F_{2m+1,2m+1}(x)&=& \sum_{i=0}^{k-1}a_i(i+x)^{2m+1}-\sum_{i=0}^{k-1}
a_i(2k-i+x)^{2m+1}\\
&=& F_{2m,2m+1}(x)+F_{2m,2m+1}(-2k-x).
\end{eqnarray*}
$F_{2m+1,2m+1}$ is a constant function, so
$$
F_{2m+1,2m+1}(x)=F_{2m+1,2m+1}(-k)=2F_{2m,2m+1}(-k).
$$
Since $a_i=a_{k-i}$, we have that
$$
F_{2m,2m+1}(-k)=-F_{2m,2m+1}(0).
$$
By induction hypothesis $F_{2m,2m+1}$ is a linear function (thus 1--1), so
$F_{2m,2m+1}(-k)\not=0$.
\item $n=2m+1$. Similarly to the previous case one can easily deduce that
$$
F_{2m+2,2m+2}(x)=2F_{2m+1,2m+2}(-(2k+1)/2)
$$
and
$$
F_{2m+1,2m+2}(-(2k+1)/2)=-F_{2m+1,2m+2}(1/2),
$$
which combined with the induction hypothesis implies that
$$
F_{2m+1,2m+2}(-(2k+1)/2) \not=0.
$$
\hfill $\square$
\end{itemize}

\begin{cor}
Let $P_n=\langle a_0,\dots,a_k\rangle$ be a $P$-sequence. Then:
\begin{itemize}
\item $\sum_{i=1}^ka_ii^s=0,\ s=0,\dots,n-1$;
\item $\sum_{i=1}^ka_ii^n\not=0$.
\end{itemize}
\end{cor}
{\bf Proof.} This is an immediate consequence of Theorem 1 and the
fact that
$$
\sum_{i=1}^ki^s=F_{n,s}(0).
$$
\hfill $\square$
\vspace{3ex}

Now let us describe how can one use the $P$-sequences in order to generate the
Prouhet-Tarry-Escott pairs:
\begin{itemize}
\item[(A)] Let $n\geqslant 2$ be an arbitrary integer and let $P_n=\langle a_0,\dots
,a_k\rangle$ be the $n$-th $P$-sequence. By Lemma 1  we have that
$$
\sum_{i=0}^ka_i(pi+l)^s=0,\ \ s=0,\dots,n-1, \ l\in\mathbb{Z},\ \
p\in\mathbb{Z}\setminus\{0\}
$$
(observe that $\sum\limits_{i=0}^ka_i(pi+l)^s=p^sF_{n,s}(l/p)$).
Since each $P$-sequence has the same number of $1$s and
$-1$s, we have that sets $X_{p,l}$ and $Y_{p,l}$ defined by
$$
X_{p,l}=\{pi+l\ |\ 0\leqslant i\leqslant k\land a_i=-1\}\ {\rm and}\
Y_{p,l}=\{pi+l\ |\ 0\leqslant i\leqslant k\land a_i=1\}
$$
have the same (finite) cardinality and
$$
\sum_{x\in X_{p,l}}x^s=\sum_{y\in Y_{p,l}}y^s,\ \ s=0,\dots,n-1,
$$
which at once yields a solution of the Prouhet-Tarry-Escott problem.
\item[(B)] For the $n$-th $P$-sequence $P_n=\langle a_0,\dots,a_k\rangle$ let
$$
X_n=\{i\in\mathbb{N}^+\ |\ i\leqslant k\land a_i=-1\}\ {\rm and}\
Y_n=\{i\in\mathbb{N}^+\ |\ i\leqslant k\land a_i=1\}.
$$
Clearly, $X_n$ and $Y_n$ are disjoint and $|X_n|=|Y_n|+1$.
Furthermore, using the definition of the notion of a $P$-sequence
one can easily check that
$$
X_{2n+1}\subset X_{2n+2}\ {\rm and}\ Y_{2n+1}\subset Y_{2n+2},
$$
and the sets $U=X_{2n+2}\setminus X_{2n+1}$ and
$V=Y_{2n+2}\setminus Y_{2n+1}$ are disjoint and have the same
cardinality. Bearing in mind the Corollary 1, we see that for each
$s\in \{0,\dots,2n\}$
$$
\sum_{i\in U}i^s=\sum_{i\in X_{2n+2}}i^s-\sum_{i\in X_{2n+1}}i^s=
\sum_{i\in Y_{2n+2}}i^s-\sum_{i\in Y_{2n+1}}i^s=\sum_{i\in V}i^s.
$$
For instance, if $n=4$, then $X_4=\{1,2,6,7,11,12\}$, $Y_4=\{4,5,8,9,13\}$,
$U=\{7,11,12\}$, $V=\{8,9,13\}$ and
$$
7^s+11^s+12^s=8^s+9^s+13^s,\ \ s=1,2.
$$
\end{itemize}
Observe that (A) and (B) are not equivalent (i.e. they produce different
sets of Prouhet-Tarry-Escott pairs).

\end{document}